\documentclass[12pt,reqno]{elsarticle}

\usepackage{amssymb}

\usepackage{natbib}
\usepackage[colorlinks=true,linkcolor=black, citecolor=blue, urlcolor=blue]{hyperref}
\usepackage{epsfig}
\usepackage{graphics}
\usepackage{amsfonts}
\usepackage{amsmath}
\usepackage{latexsym,mathrsfs}
\usepackage[nottoc, notlof, notlot]{tocbibind}
\newtheorem{theorem}{Theorem}[section]
\newtheorem{lemma}{Lemma}[section]
\newtheorem{proposition}{Proposition}[section]
\newtheorem{corollary}{Corollary}[section]

\newtheorem{definition}{Definition}[section]
\newtheorem{example}{Example}[section]

\newtheorem{remark}{Remark}[section]

\DeclareMathOperator{\infess}{ess\,inf}
\DeclareMathOperator{\supess}{ess\,sup}
\DeclareMathOperator{\supp}{supp}
\DeclareMathOperator{\Dim}{Dim}

\numberwithin{equation}{section}
\DeclareMathOperator{\R}{\mathbb{R}}
\DeclareMathOperator{\N}{\mathbb{N}}

\journal{Boletim da Sociedade Paranaense de Matem\'{a}tica}
\begin{document}
\begin{frontmatter}


\title{Multifractal dimensions for projections of measures}
\author{ Bilel
SELMI}
\address{bilel.selmi@fsm.rnu.tn}
\address{%
Faculty of Sciences of Monastir\\
Department of Mathematics\\
5000-Monastir\\
Tunisia}


\date{May 20, 2011}
\begin{abstract}
In this paper, we study the multifractal Hausdorff and packing
dimensions of Borel probability measures and study their behaviors
under orthogonal projections.  In particular, we try through these
results to improve  the main result of M. Dai in \cite{D} about the
multifractal analysis of a measure of multifractal exact dimension.

\end{abstract}
\begin{keyword}
 Multifractal analysis, Dimensions of
measures, Projection. \MSC[2010] {28A78, 28A80.}
\end{keyword}
\end{frontmatter}

\addcontentsline{toc}{section}{Multifractal dimensions for
projections of measures}
\section{Introduction}
The notion of dimensions is an important tool in the classification
of subsets in $\mathbb{R}^n$. The Hausdorff and packing dimensions
appear as some of the most common examples in the literature. The
determination of set's dimensions is naturally connected to the
auxilliary Borel measures supported by these sets. Moreover, the
estimation of a set's dimension is naturally related to the
dimension of a probability measure $\nu$ in $\mathbb{R}^n$. In this
way, thinking particularly to sets of measure zero or one, leads to
the respective definitions of the lower and upper Hausdorff
dimensions of $\nu$ as follows
 $$
\underline{\dim}(\nu)=\inf\Big\{\dim(E);\; E \subseteq
\mathbb{R}^n\; \text{and}\; \nu(E)>0\Big\}
 $$
and
$$
\overline{\dim}(\nu)=\inf\Big\{\dim(E);\; E \subseteq \mathbb{R}^n\;
\text{and}\;\nu(E)=1\Big\},
 $$
where $\dim(E)$ denotes the Hausdorff dimension of $E$
(see\cite{F}). If $\underline{\dim}(\nu)= \overline{\dim}(\nu)$,
this common value is denoted by ${\dim}(\nu)$. In this case, we say
that $\nu$ is unidimensional. Similarly, we define respectively the
lower and upper packing dimensions of $\nu$ by
 $$
\underline{\Dim}(\nu)=\inf\Big\{\Dim(E);\;E \subseteq \mathbb{R}^n\;
\text{and}\; \nu(E)>0\Big\}
 $$
and
 $$
\overline{\Dim}(\nu)=\inf\Big\{\Dim(E);\;E \subseteq \mathbb{R}^n\;
\text{and}\; \nu(E)=1\Big\},
 $$
where $\Dim(E)$ is the packing dimension of $E$ (see \cite{F}).
Also, if the equality $\underline{\Dim}(\nu)= \overline{\Dim}(\nu)$
is satisfied, we denote by ${\Dim}(\nu)$ their common value.

\bigskip
The lower and upper Hausdorff dimensions of $\nu$ were studied by
A.H. Fan in \cite{FF, FF1}. They are related to the Hausdorff
dimension of the support of $\nu$. A similar approach, concerning
the packing dimensions, was developed by Tamashiro in \cite{T}.
There are numerous works in which estimates of the dimension of a
given measure are obtained \cite{NBC, B, D, F, FLR, H1, H2, H, L,
BBSS}. When $\overline{\dim}(\nu)$ is small (resp.
$\underline{\dim}(\nu)$ is large), it means that  $\nu$ is singular
(resp. regular) with respect to the Hausdorff measure. Similar
definitions are used when
concerned with the upper and lower packing dimensions.\\
Note that, in many works (see for example \cite{F, H1, H2, H}), the
quantities $\underline{\dim}(\nu)$, $\overline{\dim}(\nu)$,
$\underline{\Dim}(\nu)$ and $\overline{\Dim}(\nu)$ are related to
the asymptotic behavior of the function $\alpha_\nu(x,r)=
\frac{\log\nu(B(x,r))}{\log r}$.

\bigskip
One of the main problems in multifractal analysis is to understand
the multifractal spectrum, the R\'{e}nyi dimensions and their
relationship with each other. During the past 20 years, there has
been enormous interest in computing the multifractal spectra of
measures in the mathematical literature and within the last 15 years
the multifractal spectra of various classes of measures in Euclidean
space $\mathbb{R}^n$ exhibiting some degree of self-similarity have
been computed rigorously (see  \cite{F, LO, Pe} and the references
therein). In an attempt to develop a general theoretical framework
for studying the multifractal structure of arbitrary measures, Olsen
\cite{LO} and Pesin \cite{Pes} suggested various ways of defining an
auxiliary measure in very general settings.  For more details and
backgrounds on multifractal analysis and its applications, the
readers may be referred also to the following essential references
\cite{NB, NBC, BJ, BB, BBH, BenM, CO, DB, BD1, FM1, MMB, MMB1, LO,
O2, O1, Ol, SH1,  SELMI1, SELMI, W, W1, W3, W4}.


 \bigskip
In this paper, we give a multifractal generalization of the results
about Hausdorff and  packing dimension of measures. We first
estimate the multifractal Hausdorff and packing dimensions of a
Borel probability measure. We try through these results to improve
the main result of M. Dai in \cite[Theorem \;A]{D} about the
multifractal analysis of a measure of exact multifractal  dimension.
We are especially based on the multifractal formalism developed by
Olsen in \cite{LO}. Then, we investigate a relationship between the
multifractal dimensions of a measure $\nu$ and its projections onto
a lower dimensional linear subspace.

\section{Preliminaries}
 We start by recalling the multifractal
formalism introduced by Olsen in \cite{LO}.  This formalism was
motivated by Olsen's wish to provide a general mathematical setting
for the ideas present in the physics literature on multifractals.

 \bigskip
Let $E\subset \mathbb{R}^n$ and $\delta>0$, we say that a collection
of balls $\big(B(x_i, r_i)\big)_i$ is a centered $\delta$-packing of
$E$ if
 $$
\forall i,\; 0<r_i<\delta,\quad x_i\in E,\; \text{and} \quad
B(x_i,r_i)\cap B(x_j, r_j)=\emptyset,\quad{\forall}\; i\neq j.
 $$
Similarly, we say that $\big(B(x_i, r_i)\big)_i$ is a centered
$\delta$-covering of $E$ if
 $$
\forall i,\; 0<r_i<\delta, \quad x_i\in E, \quad \text{and} \qquad
E\subset \bigcup_i \; B(x_i, r_i).
 $$

Let $\mu$ be a Borel probability measure on $\mathbb{R}^n$. For $q,
t\in\mathbb{R}$, $E \subseteq{\mathbb R}^n$ and $\delta>0$, we
define
 $$
\overline{{\mathcal P}}^{q,t}_{\mu,\delta}(E) =\displaystyle \sup
\left\{\sum_i \mu(B(x_i,r_i))^q (2r_i)^t\right\},
 $$
where the supremum is taken over all centered $\delta$-packings of
$E$. The generalized packing pre-measure is  given by
 $$
\overline{{\mathcal P}}^{q,t}_{\mu}(E)
=\displaystyle\inf_{\delta>0}\overline{{\mathcal
P}}^{q,t}_{\mu,\delta}(E).
 $$
In a similar way, we define
 $$
\overline{{\mathcal H}}^{q,t}_{\mu,\delta}(E) = \displaystyle\inf
\left\{\sum_i \mu(B(x_i,r_i))^q(2r_i)^t\right\},
 $$
where the infinimum is taken over all centered $\delta$-coverings of
$E$. The generalized Hausdorff pre-measure is defined by
 $$
\overline{{\mathcal H}}^{q,t}_{\mu}(E) = \displaystyle\sup_{
\delta>0}\overline{{\mathcal H}}^{q,t}_{\mu,\delta}(E).
 $$
Especially, we have the conventions $0^q=\infty$ for $q\leq0$ and
$0^q=0$ for $q>0$.

 \bigskip
Olsen \cite{LO} introduced the following modifications on the
generalized Hausdorff and packing measures,
$$
{\mathcal H}^{q,t}_{\mu}(E)=\displaystyle\sup_{F\subseteq
E}\overline{{\mathcal H}}^{q,t}_{\mu}(F)\quad\text{and}\quad
{\mathcal P}^{q,t}_{\mu}(E) = \inf_{E \subseteq \bigcup_{i}E_i}
\sum_i \overline{\mathcal P}^{q,t}_{\mu}(E_i).
 $$

The functions ${\mathcal H}^{q,t}_{\mu}$ and ${\mathcal
P}^{q,t}_{\mu}$ are metric outer measures and thus measures on the
family of Borel subsets of $\mathbb{R}^n$. An important feature of
the Hausdorff and packing measures is that ${\mathcal
P}^{q,t}_{\mu}\leq{\overline{\mathcal P}}^{q,t}_{\mu}$. Moreover,
there exists an integer $\xi\in\mathbb{N}$, such that ${\mathcal
H}^{q,t}_{\mu}\leq\xi{\mathcal P}^{q,t}_{\mu}.$ The measure
${\mathcal H}^{q,t}_{\mu}$ is of course a multifractal
generalization of the centered Hausdorff measure, whereas ${\mathcal
P}^{q,t}_{\mu}$ is a multifractal generalization of the packing
measure. In fact, it is easily seen that, for  $t\geq0$,  one has
$$2^{-t} {\mathcal H}^{0,t}_{\mu}\leq {\mathcal H}^{t}\leq {\mathcal
H}^{0,t}_{\mu}\quad\text{ and}\quad{\mathcal
P}^{0,t}_{\mu}={\mathcal P}^{t},$$ where ${\mathcal H}^{t}$ and
${\mathcal P}^{t}$ denote respectively the $t$-dimensional Hausdorff
and  $t$-dimensional packing measures.

 \bigskip
The measures ${\mathcal H}^{q,t}_{\mu}$ and ${\mathcal
P}^{q,t}_{\mu}$ and the pre-measure ${\overline{\mathcal
P}}^{q,t}_{\mu}$ assign in a usual way a multifractal dimension to
each subset $E$ of $\mathbb{R}^n$. They are respectively denoted by
$\dim_{\mu}^q(E)$, $\Dim_{\mu}^q(E)$ and $\Delta_{\mu}^q(E)$ (see
\cite{LO}) and satisfy
 $$
\begin{array}{lllcr}
\dim_{\mu}^q(E) &=&\inf \Big\{ t\in\R; \quad {\mathcal
H}^{{q},t}_{\mu}(E) =0\Big\}=\sup \Big\{ t\in\R; \quad {\mathcal
H}^{{q},t}_{\mu}(E) =+\infty\Big\}, \\ \\
\Dim_{\mu}^q(E) &=& \inf \Big\{  t\in\R; \quad {\mathcal
P}^{{q},t}_{\mu}(E) =0\Big\}=\sup \Big\{  t\in\R; \quad {\mathcal
P}^{{q},t}_{\mu}(E) =+\infty\Big\}, \\ \\
\Delta_{\mu}^q(E) &=& \inf \Big\{ t\in\R; \quad \overline{\mathcal
P}^{{q},t}_{\mu}(E) =0\Big\}=\sup\Big\{ t\in\R; \quad
\overline{\mathcal P}^{{q},t}_{\mu}(E) =+\infty\Big\}.
\end{array}
 $$
The number $\dim_{\mu}^q(E)$ is an obvious multifractal analogue of
the Hausdorff dimension $\dim(E)$ of $E$ whereas $\Dim_{\mu}^q(E)$
and $\Delta_{\mu}^q(E)$ are obvious multifractal analogues of the
packing dimension $\Dim(E)$ and the pre-packing dimension
$\Delta(E)$ of $E$ respectively. In fact, it follows immediately
from the definitions that
$$
\dim(E)=\dim_{\mu}^0(E),\;\;\;\Dim(E)=\Dim_{\mu}^0(E)\quad\text{and}\quad\Delta(E)=\Delta_{\mu}^0(E).
$$

\section{Multifractal Hausdorff and packing dimensions of measures}
Now, we introduce the multifractal analogous of the Hausdorff and
packing dimensions of a Borel probability measure.

\begin{definition} The lower and upper multifractal Hausdorff
dimensions of a measure $\nu$ with respect to a measure $\mu$ are
defined by
 $$
\underline{\dim}_{\mu}^q(\nu)=\inf\Big\{\dim_{\mu}^q(E);\; E
\subseteq \mathbb{R}^n\; \text{and}\; \nu(E)>0\Big\}
 $$
and
 $$
\overline{\dim}_{\mu}^q(\nu)=\inf\Big\{\dim_{\mu}^q(E);\; E
\subseteq \mathbb{R}^n\; \text{and}\; \nu(E)=1\Big\}.
 $$
We denote by ${\dim}_{\mu}^q(\nu)$ their common value, if the
equality $\underline{\dim}_{\mu}^q(\nu)= \overline{\dim}_{\mu}^q
(\nu)$ is satisfied.
\end{definition}

\begin{definition} The lower and upper multifractal packing
dimensions of a measure $\nu$ with respect to a measure $\mu$ are
defined by
 $$
\underline{\Dim}_{\mu}^q(\nu)=\inf\Big\{\Dim_{\mu}^q(E);\; E
\subseteq \mathbb{R}^n\; \text{and}\; \nu(E)>0\Big\}
 $$
and
$$
\overline{\Dim}_{\mu}^q(\nu)=\inf\Big\{\Dim_{\mu}^q(E);\; E
\subseteq \mathbb{R}^n\; \text{and}\; \nu(E)=1\Big\}.
 $$
When $\underline{\Dim}_{\mu}^q(\nu)= \overline{\Dim}_{\mu}^q(\nu)$,
we denote by ${\Dim}_{\mu}^q(\nu)$ their common value.
\end{definition}

\begin{definition}Let $\mu, \nu$ be two Borel probability measures on $\mathbb{R}^n$.\noindent
\begin{enumerate}
\item
We say that $\mu$ is absolutely continuous with respect to $\nu$ and
write $\mu\ll \nu$ if, for any set $A\subset\R^n$,
$\nu(A)=0\Rightarrow \mu(A)=0$.

\item $\mu$ and $\nu$ are said to be mutually singular and we write
$\mu\bot \;\nu$ if there exists a set $A\subset\R^n$, such that
$\mu(A)=0=\nu(\R^n\setminus A).$
\end{enumerate}
\end{definition}

The quantities $\underline{\dim}_\mu^q(\nu)$ and
$\overline{\dim}_\mu^q(\nu)$ \big(resp. $\underline{\Dim}_\mu^q
(\nu)$ and $\overline{\Dim}_\mu^q(\nu)$\big) allow to compare the
measure $\nu$ with the generalized Hausdorff (resp. packing)
measure. More precisely, we have the following result.
\begin{theorem}\label{0} Let $\mu, \nu$ be two Borel probability measures on $\mathbb{R}^n$ and
$q\in \mathbb{R}$. We have,

\begin{enumerate}
\item $\underline{\dim}_{\mu}^q(\nu)=\sup\Big\{t\in\R;\;\nu\ll{\mathcal
H}^{q,t}_{\mu}\Big\}$\; and\;
$\overline{\dim}_{\mu}^q(\nu)=\inf\Big\{t\in\R;\; \nu\bot{\mathcal
H}^{q,t}_{\mu}\Big\}.$

  \bigskip
\item$ \underline{\Dim}_{\mu}^q(\nu)=\sup\Big\{t\in\R;\; \nu\ll{\mathcal
P}^{q,t}_{\mu}\Big\}$ \;and\; $
\overline{\Dim}_{\mu}^q(\nu)=\inf\Big\{t\in\R;\; \nu\bot{\mathcal
P}^{q,t}_{\mu}\Big\}.$
\end{enumerate}
\end{theorem}
\noindent{\bf Proof.}
\begin{enumerate}
\item Let's prove that $\underline{\dim}_{\mu}^q(\nu) =\sup\Big\{
t\in\R;\;\nu\ll{\mathcal H}^{q,t}_{\mu}\Big\}$.
Define
 $$
s=\sup\Big\{t\in\R;\; \nu\ll{\mathcal H}^{q,t}_{\mu}\Big\}.
 $$

For any $t<s$ and $E\subseteq\R^n$, such that $\nu(E)>0$, we have
${\mathcal H}^{q,t}_{\mu}(E)>0$. It follows that
$\dim_{\mu}^q(E)\geq t$ and then, $\underline{\dim}_{\mu}^q(\nu)\geq
t$. We deduce that $\underline{\dim}_{\mu}^q(\nu)\geq s$.\\

On the other hand, for any $t>s$, there exists a set
$E\subseteq\R^n$, such that $\nu(E)>0$ and ${\mathcal
H}^{q,t}_{\mu}(E)=0$. Consequently, $\dim_{\mu}^q(E)\leq t$ and so,
$\underline{\dim}_{\mu}^q(\nu)\leq t$. This leads to
$\underline{\dim}_{\mu}^q(\nu)\leq s$.

  \bigskip
Now, we prove that $\overline{\dim}_{\mu}^q(\nu)=\inf\Big\{t\in\R;\;
\nu\bot{\mathcal H}^{q,t}_{\mu}\Big\}.$ For this, we define
 $$
s'=\inf\Big\{t\in\R;\; \nu\bot{\mathcal H}^{q,t}_{\mu}\Big\}.
 $$
For $t>s'$, there exists a set $E\subseteq\R^n$, such that
${\mathcal H}^{q,t}_{\mu}(E)=0=\nu(\R^n\setminus E)$. Then,
$\dim_{\mu}^q(E)\leq t$. Since $\nu(E)=1$, then
$\underline{\dim}_{\mu}^q(\nu)\leq t$
and $\underline{\dim}_{\mu}^q(\nu)\leq s'$.\\

Now, for $t<s'$, take $E\subseteq\R^n$, such that ${\mathcal
H}^{q,t}_{\mu}(E)>0$ and $\nu(E)=1$. It can immediately seen that
$\dim_{\mu}^q(E)\geq t$. Then, $\underline{\dim}_{\mu}^q(\nu)\geq
t$. It follows that $\underline{\dim}_{\mu}^q(\nu)\geq s'$. This
ends the proof of assertion (1).

  \bigskip
\item The proof of assertion (2) is given in \cite[Theorem 2]{L}.$\hfill\square$
\end{enumerate}
\begin{remark}
When the upper multifractal Hausdorff (resp. packing) dimension of
the measure is small, it means that the measure $\nu$ is ``very
singular" with respect to the generalized multifractal Hausdorff
(resp. packing) measure. In the same way, when the lower
multifractal (resp. packing) dimension of the measure is large, then
the measure $\nu$ is ``quite regular" with respect to the
generalized multifractal Hausdorff (resp. packing) measure.
\end{remark}

 \bigskip \bigskip
The quantities $\underline{\dim}_{\mu}^q(\nu)$,
$\overline{\dim}_{\mu}^q(\nu)$, $\underline{\Dim}_{\mu}^q(\nu)$ and
$\overline{\Dim}_{\mu}^q(\nu)$ are related to the asymptotic
behavior of the function ${\alpha}_{\mu,\nu}^q(x,r)$, where
 $$
{\alpha}_{\mu,\nu}^q(x,r)=\frac{\log \nu\big(B(x,r)\big) -q\log\mu\big(B(x,r)\big) }{\log r}.
 $$

Notice that the chararcterization of the lower und upper packing
dimensions by the function $\alpha_{\mu,\nu}^q$ is proved by J. Li
in \cite[Theorem 3]{L}. In the following theorem we prove similar
results for the Hausdorff dimensions.
\begin{theorem}\label{1} Let $\mu, \nu$ be two Borel probability measures on $\mathbb{R}^n$ and
$q\in \mathbb{R}$. Let
$$
\underline{\alpha}_{\mu,\nu}^q(x)=\displaystyle\liminf_{ r\to 0}\;
{\alpha}_{\mu,\nu}^q(x,r)\quad\text{and}\quad
\overline{\alpha}_{\mu,\nu}^q(x)=\displaystyle\limsup_{ r\to 0}\;
{\alpha}_{\mu,\nu}^q(x,r).
 $$
We have,

\begin{enumerate}
\item
 $
\underline{\dim}_{\mu}^q(\nu)=\infess\underline{\alpha}_{\mu,\nu}^q(x)
\quad\text{and}\quad
\overline{\dim}_{\mu}^q(\nu)=\supess\underline{\alpha}_{\mu,\nu}^q(x).
 $

  \bigskip
\item
$
\underline{\Dim}_{\mu}^q(\nu)=\infess\overline{\alpha}_{\mu,\nu}^q(x)\quad\text{and}\quad
\overline{\Dim}_{\mu}^q(\nu)=\supess\overline{\alpha}_{\mu,\nu}^q(x),
 $
\end{enumerate}
where the essential bounds being related to the measure $\nu$.
\end{theorem}

\noindent{\bf Proof.}
\begin{enumerate}
\item We prove that $\underline{\dim}_{\mu}^q(\nu)
=\text{ess}\inf\underline{\alpha}_{\mu,\nu}^q(x)$.

 \bigskip
Let $\alpha<\text{ess}\inf\underline{\alpha}_{\mu,\nu}^q(x)$. For
$\nu$-almost every $x$, there exists $r_0>0$, such that $0<r<r_0$
and
 $$
\nu(B(x,r))<\mu(B(x,r))^q ~~r^\alpha.
 $$
Denote by
 $$
F_n=\left\{x;\;\nu(B(x,r))<\mu(B(x,r))^q ~~r^\alpha,\;\,
\text{for}\;\,0<r<\frac1n\right\}.
 $$
Let $F=\cup_n F_n$. It is clear that $\nu(F)=1$. Take $E$ be a Borel
subset of $\R^n$ satisfying $\nu(E)> 0$. We have $\nu(E\cap F)>0$
and there exists an integer $n$, such that $\nu(E\cap F_n)>0$.

Let $\delta>0$ and $\big(B(x_i,r_i)\big)_i$ be a centered
$\delta$-covering of $E \cap F_n$. We have
 $$
\displaystyle\sum_{i}\nu(B(x_{i},r_{i}))\leq 2^{-\alpha}
\displaystyle\sum_{i}\mu(B(x_{i},r_{i}))^q (2r_i)^\alpha,
 $$
so that
 $$
2^{\alpha}\nu(E\cap F_n)\leq \overline{{\mathcal
H}}^{q,\alpha}_{\mu,\delta}(E\cap F_n).
 $$
Letting $\delta \to 0$ gives that
 $$
2^{\alpha}\nu(E\cap F_n)\leq \overline{{\mathcal
H}}^{q,\alpha}_{\mu}(E\cap F_n)\leq {\mathcal
H}^{q,\alpha}_{\mu}(E\cap F_n).
 $$
It follows that
 $$
{\mathcal H}^{q,\alpha}_{\mu}(E)\geq{\mathcal
H}^{q,\alpha}_{\mu}(E\cap F_n)>0\;\Rightarrow\; \dim_{\mu}^q(E)\geq
\alpha.
 $$
We have proved that
$$\underline{\dim}_{\mu}^q(\nu)\geq\text{ess}\inf\underline{\alpha}_{\mu,\nu}^q(x).$$

  \bigskip
On the other hand, if
$\text{ess}\inf\underline{\alpha}_{\mu,\nu}^q(x)=\alpha$. For
$\varepsilon>0$, let
 $$
E_\varepsilon=\Big\{x\in\supp\nu;\;
\underline{\alpha}_{\mu,\nu}^q(x)<\alpha+\varepsilon\Big\}.
 $$
It is clear that $\nu(E_\varepsilon)>0$. This means that
$\underline{\dim}_{\mu}^q(\nu)\leq{\dim}_{\mu}^q(E_\varepsilon)$. We
will prove that
 $$
{\dim}_{\mu}^q(E_\varepsilon)\leq \alpha+\varepsilon,  \;\;
\forall\; \varepsilon>0.
 $$
Let $E\subset E_\varepsilon$ and $x\in E$. Then, for all $\delta>0$
we can find $0<r_x<\delta$, such that
 $$
\nu(B(x,r_x))>\mu(B(x,r_x))^q ~~r_x^{\alpha+\varepsilon}.
 $$
Take $\delta>0$. The family $\Big(B(x,r_x)\Big)_{x\in E}$ is a
centered $\delta$-covering of $ E.$ Using Besicovitch's Covering
Theorem (see \cite{F,M1}), we can construct $\xi$ finite or
countable sub-families
$\Big(B(x_{1j},r_{1j})\Big)_j$,....,$\Big(B(x_{\xi j},r_{\xi
j})\Big)_j$, such that each  $E$ satisfies
 $$
E\subseteq\displaystyle\bigcup_{i=1}^\xi\bigcup_jB(x_{ij},r_{ij})
\quad \text{and} \quad \Big(B(x_{ij},r_{ij})\Big)_j\quad \text{is a
} \delta\text{-packing of }E.
 $$
We get
\begin{eqnarray*}
\displaystyle\sum_{i,j}\mu(B(x_{ij},r_{ij}))^q
(2r_{ij})^{\alpha+\varepsilon} &\leq& \xi 2^{\alpha+\varepsilon}
\sum_j \nu(B(x_{ij},r_{ij}))\leq \xi 2^{\alpha+\varepsilon}
\nu(\mathbb{R}^n).
\end{eqnarray*}
Consequently,
 $$
\overline{\mathcal H}^{q,\alpha+\varepsilon}_{\mu,\delta}(E)\leq\xi
2^{\alpha+\varepsilon} \nu(\mathbb{R}^n)\;\Rightarrow\;
\overline{\mathcal H}^{q,\alpha+\varepsilon}_{\mu}(E)\leq\xi
2^{\alpha+\varepsilon} \nu(\mathbb{R}^n).
 $$
We obtain thus
 $$
{\mathcal H}^{q,\alpha+\varepsilon}_{\mu}(E_\varepsilon)\leq\xi
2^{\alpha+\varepsilon} \nu(\mathbb{R}^n)<\infty.
 $$
Therefore,
 $$
\dim_{\mu}^q(E)\leq \alpha+\varepsilon\quad\text{and}\quad
\underline{\dim}_{\mu}^q(\nu)\leq\text{ess}\inf\underline{\alpha}_{\mu,\nu}^q(x).
 $$
We prove in a similar way that
$\overline{\dim}_{\mu}^q(\nu)=\text{ess}\sup\underline{\alpha}_{\mu,\nu}^q(x)$.
\end{enumerate}

\begin{corollary}\label{cor1}
Let $\mu, \nu$ be two Borel probability measures on $\mathbb{R}^n$
and take $q, \alpha\in \mathbb{R}$. We have,

\begin{enumerate}
\item $\underline{\dim}_{\mu}^q(\nu)\geq\alpha$ \;if and only if\;
$\underline{\alpha}_{\mu,\nu}^q(x)\geq\alpha$\; for\; $\nu$-a.e.
$x$.

 \smallskip
\item $\overline{\dim}_{\mu}^q(\nu)\leq\alpha$ \;if and only if\;
$\underline{\alpha}_{\mu,\nu}^q(x)\leq\alpha$ \;for\; $\nu$-a.e.
$x$.

 \smallskip
\item $\underline{\Dim}_{\mu}^q(\nu)\geq\alpha$ \;if and only if\;
$\overline{\alpha}_{\mu,\nu}^q(x)\geq\alpha$ for $\nu$-a.e. $x$.

 \smallskip
\item $\overline{\Dim}_{\mu}^q(\nu)\leq\alpha$ \;if and only if\;
$\overline{\alpha}_{\mu,\nu}^q(x)\leq\alpha$ for $\nu$-a.e. $x$.
\end{enumerate}
\end{corollary}

\noindent{\bf Proof.} Follows immediately from Theorem \ref{1}.

\begin{example}
\end{example} We recall the definition of the deranged Cantor set
(see \cite{B11, B1, B12, BB}).

\noindent Let $I_\emptyset= [0, 1]$. We obtain respectively the left
and right sub-intervals $I_{\varepsilon,1}$ and $I_{\varepsilon,2}$
of $I_\varepsilon$ by deleting the middle open sub-interval of
$I_\varepsilon$ inductively for each $\varepsilon\in\{1, 2\}^n$,
where $n \in\N$.

\noindent We consider the sequence
$$
\mathcal{C}_n =\displaystyle\bigcup_{\varepsilon\in\{1, 2\}^n}
I_\varepsilon.
 $$
$\{\mathcal{C}_n \}_{n\in\N}$ is a decreasing sequence of closed
sets.

 \smallskip \smallskip
\noindent For each $n\in \N$ and each $\varepsilon\in\{1, 2\}^n$, we
put
 $$
|I_{\varepsilon,1} | / |I_{\varepsilon} |=
c_{\varepsilon,1}\quad\text{and}\quad |I_{\varepsilon,2} | /
|I_{\varepsilon} |= c_{\varepsilon,2},
 $$
where $|I|$ is the diameter of $I$. The set
$\mathcal{C}=\displaystyle\bigcap_{n\geq0}\mathcal{C}_n$ is called a
deranged Cantor set.

 \bigskip
Let $\nu$ be a probability measure supported by the deranged Cantor
set $\mathcal{C}$ and $\mu$ be the Lebesgue measure on
$I_\emptyset$. For $\varepsilon_1,...,\varepsilon_n\in\{1,2\}$, we
denote by $I_{\varepsilon_1,...,\varepsilon_n}$ the basic set of
level $n$. For $x\in C$, we denote by $I_n(x)$ the $n$-th level set
containing $x$. We introduce the sequence of random variables $X_n$
defined by
 $$
X_n(x)=-\log_3\left(\frac{\nu(I_n(x))}{\nu(I_{n-1}(x))}\right).
 $$
We have
 $$
\frac{S_n(x)}{n}=\frac{X_1(x)+...+X_n(x)}{n}=
\frac{\log(\nu(I_n(x))}{\log\mid I_n(x)\mid}.
 $$
By Lemma 1 in \cite{BB}, we have for all $x\in \mathcal{C}$,
 $$
\liminf_{n\to\infty}\frac{\log(\nu(I_n(x))}{\log\mid
I_n(x)\mid}=\liminf_{r\to0}\frac{\log(\nu(B(x,r))}{\log r}
 $$
and
 $$
\limsup_{n\to\infty}\frac{\log(\nu(I_n(x))}{\log\mid
I_n(x)\mid}=\limsup_{r\to0}\frac{\log(\nu(B(x,r))}{\log r}.
 $$
The quantities $\underline{\dim}_{\mu}^q(\nu)$ and
$\overline{\dim}_{\mu}^q(\nu)$ are related to the asymptotic
behavior of the sequence $\displaystyle\frac{S_n}{n}$. More
precisely, we have the following two relations
 $$
\underline{\dim}_{\mu}^q(\nu)
=\text{ess}\inf\left\{\liminf_{n\to\infty}
\frac{S_n(x)}{n}-q\right\} \quad\text{and}\quad
\overline{\dim}_{\mu}^q(\nu)
=\text{ess}\sup\left\{\liminf_{n\to\infty}
\frac{S_n(x)}{n}-q\right\}.
 $$
In the same way, we can also prove that
 $$
\underline{\Dim}_{\mu}^q(\nu)
=\text{ess}\inf\left\{\limsup_{n\to\infty}\frac{S_n(x)}{n}-q\right\}
\quad\text{and}\quad \overline{\Dim}_{\mu}^q(\nu)
=\text{ess}\sup\left\{\limsup_{n\to\infty}\frac{S_n(x)}{n}-q\right\}.
 $$

We say that the measure $\nu$ is $(q, \mu)$-unidimensional if
$\overline{\dim}_{\mu}^q(\nu)=\underline{\dim}_{\mu}^q(\nu)$. We
also say that $\nu$ has an exact multifractal packing dimension
whenever $\overline{\Dim}_{\mu}^q(\nu)
=\underline{\Dim}_{\mu}^q(\nu)$. In general, a Borel probability
measure is not $(q,\mu)$-unidimensional and
$\overline{\Dim}_{\mu}^q(\nu)\neq\underline{\Dim}_{\mu}^q(\nu)$.

 \bigskip
In the following, we are interested to the $(q,
\mu)$-unidimensionality and ergodicity of $\nu$ and to the calculus
of its multifractal Hausdorff and packing dimensions. Our purpose in
the following theorem is to prove the main Theorem of M. Dai
\cite[Theorem \;A]{D} under less restrictive hypotheses.
\begin{theorem}\label{p}
The measure $\nu$ is $(q, \mu)$-unidimensional with
${\dim}_{\mu}^q(\nu)=\alpha$ if and only if the following two
conditions are satisfied.
\begin{enumerate}
\item There exists a set $E$ of $\R^n$ with ${\dim}_{\mu}^q(E)=\alpha$, such that $\nu(E)=1$.

\item $\nu(E)=0$, for every Borel set $E$ satisfying
${\dim}_{\mu}^q(E)<\alpha$.
\end{enumerate}
\end{theorem}

\noindent{\bf Proof.} We can deduce from Theorems \ref{0} and
\ref{1} that $\nu$ is $(q, \mu)$-unidimensional if and only if we
have the following assertions.
\begin{enumerate}
\item $\nu$ is absolutely continuous with respect to ${\mathcal
H}^{q,\alpha-\varepsilon}_{\mu}$, for all $\varepsilon>0$.

\item $\nu$ and ${\mathcal H}^{q,\alpha+\varepsilon}_{\mu}$ are
mutually singular, for all $\varepsilon>0$.
\end{enumerate}
Then, the proof of Theorem \ref{p} becomes an easy consequence of
the following lemma.

\begin{lemma}\cite{D} The following conditions are equivalent.
\begin{enumerate}
\item We have,
\begin{enumerate}
\item there exists a set $E$ of $\R^n$ with ${\dim}_{\mu}^q(E)=\alpha$, such that $\nu(E)=1$.

\item $\nu(E)=0$, for every Borel set $E$ satisfying
${\dim}_{\mu}^q(E)<\alpha$.
\end{enumerate}

\item We have,
\begin{enumerate}
\item $\nu\ll{\mathcal
H}^{q,\alpha-\varepsilon}_{\mu}$ for all $\varepsilon>0$.
\item $\nu\bot{\mathcal
H}^{q,\alpha+\varepsilon}_{\mu}$ for all $\varepsilon>0$.
\end{enumerate}
\end{enumerate}
\end{lemma}
\begin{remark}
Theorem \ref{p} improves Dai's result \cite[Theorem \;A]{D} (we need
not to assume that $\mu$ is a doubling measure).
\end{remark}

The symmetrical results are true as well.
\begin{theorem}\label{20}
Let $\mu, \nu$ be two Borel probability measures on $\mathbb{R}^n$
and take $\alpha, q\in \R$. The following conditions are equivalent.
\begin{enumerate}
\item $\overline{\Dim}_{\mu}^q(\nu)=\underline{\Dim}_{\mu}^q(\nu)=\alpha$.
\item We have,
\begin{enumerate}
\item there exist a set $E\subset\R^n$ with
$\Dim_{\mu}^q(E)=\alpha$, such that $\nu(E)=1$,
\item if $\; E\subset\R^n$ satisfies $\Dim_{\mu}^q(E)<\alpha$, then
$\nu(E)=0$.
\end{enumerate}
\item We have,
\begin{enumerate}
\item $\nu \ll {\mathcal P}^{q,\alpha-\epsilon}_{\mu}$, for
all $\epsilon>0$.
\item $\nu \bot {\mathcal P}^{q,\alpha+\epsilon}_{\mu}$, for
all $\epsilon>0$.
\end{enumerate}
\end{enumerate}
\end{theorem}
\noindent{\bf Proof.} We can deduce from Theorems \ref{0} and
\ref{1} that the assertions (1) and (3) are equivalent. We only need
to prove the equivalence of the assertions (2) and (3).

 \bigskip
Assume that the measure $\nu$ satisfies the hypothesis (a) and (b)
of (2). Let $\; E\subset\R^n$ and suppose that ${\mathcal
P}^{q,\alpha-\epsilon}_{\mu}(E)=0$, for all $\epsilon>0$. Then, we
have that $\Dim_{\mu}^q(E)\leq \alpha-\epsilon<\alpha$. By condition
(b) of (2), we obtain $\nu(E)=0$. Thus,
 $$
\nu \ll {\mathcal P}^{q,\alpha-\epsilon}_{\mu},\quad\text{for all }
\epsilon>0.
 $$

Thanks to condition (a) of (2), there exists a set $E\subset\R^n$ of
multifractal packing dimension $\alpha$, such that $\nu(E)=1$ and
$\Dim_{\mu}^q(E)=\alpha< \alpha+\epsilon$, for all $\epsilon>0$.
Then, ${\mathcal P}^{q,\alpha+\epsilon}_{\mu}(E)=0$. Thus,
 $$
\nu \bot {\mathcal P}^{q,\alpha+\epsilon}_{\mu},\quad\text{for all }
\epsilon>0.
 $$

 \bigskip
Now, assume that $\nu$ satisfies conditions (a) and (b) of (3). This
means that $\nu \ll {\mathcal P}^{q,\alpha -\epsilon}_{\mu}$, for
all $\epsilon>0$. Taking a Borel set $E$ with $\Dim_{\mu}^q(E)
=\beta<\alpha$ and $\epsilon=\frac{\alpha-\beta}2$, we get
${\mathcal P}^{q,(\alpha+\beta)/2}_{\mu}(E)=0.$ Then, $\nu(E)=0$.

Since $\nu \bot {\mathcal P}^{q,\alpha +\epsilon }_{\mu}$, for all
$\epsilon>0$, there exists a set $F_\epsilon$ with ${\mathcal
P}^{q,\alpha +\epsilon }_{\mu}(F_\epsilon)=0$ and
$\nu(F_\epsilon)=1$. Hence,
$\Dim_{\mu}^q(F_\epsilon)\leq\alpha+\epsilon$. Choose a sequence
$(\epsilon_k)_k$ such that $\epsilon_k\to 0$ as $k\to +\infty$ and
consider the set $F=\displaystyle\bigcap_{k\geq 1} F_{\epsilon_k}$.
It is clear that $\nu(F)=1$ and
 $$
\Dim_{\mu}^q(F) \leq\displaystyle\liminf_{k\to\infty}
\Dim_{\mu}^q(F_{\epsilon_k})\leq
\displaystyle\liminf_{k\to\infty}(\alpha+\epsilon_k)=\alpha.
 $$
If $\Dim_{\mu}^q(F) =\alpha$, then the condition (a) of (2) is
satisfied for $E=F$.\\
If $\Dim_{\mu}^q(F) <\alpha$, then putting $E=F\cup G$, for some
Borel set $G$ of multifractal packing dimension $\alpha$, we obtain
 $$
\nu(E)=1 \qquad \text{and}\qquad \Dim_{\mu}^q(E) =\max \Big\{
\Dim_{\mu}^q(F), \Dim_{\mu}^q(G)\Big\} =\alpha.
 $$
\begin{proposition}
Let $\mu$  be the Lebesgue measure on $\mathbb{R}^n$, $\nu$ be a
compactly supported Borel probability measure on $\mathbb{R}^n$ and
$T:\supp\nu\rightarrow\supp\nu$ a $K$-lipschitz function. Suppose
that $\nu$ is $T$-invariant and ergodic on $\supp\nu$. Then,
 $$
\overline{\dim}_{\mu}^q(\nu)=\underline{\dim}_{\mu}^q(\nu)\quad\text{and}\quad
\overline{\Dim}_{\mu}^q(\nu)=\underline{\Dim}_{\mu}^q(\nu).
 $$
\end{proposition}
\noindent{\bf Proof.} $T$ is a $K$-lipschitz function, then
$T(B(x,r))\subseteq B(T(x), Kr)$. Since $\nu$ is $T$-invariant, then
we can deduce that
 $$
\nu\big(B(x,r)\big)\leq\nu\big(T^{-1}\big(T(B(x,r))\big)\big)\leq\nu\big(T^{-1}\big(B(T(x),
Kr)\big)\big)=\nu\big(B(T(x), Kr)\big).
 $$
It follows that,
\begin{eqnarray*}
\frac{\log \nu\big(B(x,r)\big)-q\log\mu\big(B(x,r)\big) }{\log r}
&=& \frac{\log \nu\big(B(x,r)\big) }{\log r}-q \\
&\geq& \frac{\log \nu\big(B(T(x),K r)\big) }{\log  (K r)} \times
\frac{\log  (K r)}{\log r}-q,
\end{eqnarray*}
which proves that
 $$
\underline{\alpha}_{\mu,\nu}^q(x)\geq\underline{\alpha}_{\mu,\nu}^q(T(x))\quad\text{and}
\quad\overline{\alpha}_{\mu,\nu}^q(x)\geq\overline{\alpha}_{\mu,\nu}^q(T(x)).
 $$
Since $\nu$ is ergodic, then the function
$\underline{\alpha}_{\mu,\nu}^q(x)
-\underline{\alpha}_{\mu,\nu}^q(T(x))$ \Big(resp.
$\overline{\alpha}_{\mu,\nu}^q(x)-\overline{\alpha}_{\mu,\nu}^q(T(x))$\Big)
is positive and satisfies
 $$
\int\Big(\underline{\alpha}_{\mu,\nu}^q(x)-\underline{\alpha}_{\mu,\nu}^q(T(x))\Big)
d\nu(x)=0\quad
\left(\text{resp.}\;\;\int\Big(\overline{\alpha}_{\mu,\nu}^q(x)-\overline{\alpha}_{\mu,\nu}^q(T(x))\Big)
d\nu(x)=0\right).
 $$
We can conclude that,
 $$
\underline{\alpha}_{\mu,\nu}^q(x)=\underline{\alpha}_{\mu,\nu}^q(T(x))\quad\text{and}
\quad\overline{\alpha}_{\mu,\nu}^q(x)=\overline{\alpha}_{\mu,\nu}^q(T(x))\quad
\text{for }\; \nu\; \text{-a.e. } x
 $$
and that the functions $\underline{\alpha}_{\mu,\nu}^q$,
$\overline{\alpha}_{\mu,\nu}^q$ are $T$-invariant. On the other
hand, the measure $\nu$ is ergodic and
 $$
-q\leq\underline{\alpha}_{\mu,\nu}^q\leq
\overline{\alpha}_{\mu,\nu}^q\leq n-q.
 $$
It follows that $\underline{\alpha}_{\mu,\nu}^q$
\big($\overline{\alpha}_{\mu,\nu}^q$\big) is $\nu$-almost every
where constant, which says that
 $$
\overline{\dim}_{\mu}^q(\nu)=\underline{\dim}_{\mu}^q(\nu)\quad\text{and}\quad
\overline{\Dim}_{\mu}^q(\nu)=\underline{\Dim}_{\mu}^q(\nu).
 $$
\begin{remark}
In the case where
$\overline{\alpha}_{\mu,\nu}^q(x)=\underline{\alpha}_{\mu,\nu}^q(x)=\alpha$
for $\nu$-almost all $x$, we have
${\dim}_{\mu}^q(\nu)={\Dim}_{\mu}^q(\nu)=\alpha. $  The results
developed by Heurteaux in \cite{H1, H2, H} and Fan et al. in
\cite{FF, FF1, FLR} are obtained as a special case of the
multifractal Theorems when $q$ equals $0$.
\end{remark}
\begin{example}
\end{example}
We say that the probability measure $\mu$ is a quasi-Bernoulli
measure on the Cantor set $\mathcal{C}= \{0, 1\}^{\N^*}$ if we can
find $C \geq 1$ such that
 $$
\forall x,y \in \mathcal{F}\qquad C^{-1} \mu(x)\mu(y)\leq\mu(x
y)\leq C \mu(x)\mu(y),
 $$
where $\mathcal{F}$ is the set of words written with the alphabet
$\{0, 1\}$. Let $\mathcal{F}_n$ be the set of words of length $n$,
and take $x = x_1 x_2 ... \in \mathcal{C}$, let $I_n(x)$ be the
unique cylinder $\mathcal{F}_n$ that contains $x$. Let us introduce
the function $\tau_\mu$ defined for $p\in\R$, by
 $$
\tau_\mu(p)=\displaystyle\lim_{n\to\infty}
\tau_\mu(p,n)\quad\text{with}\quad
\tau_\mu(p,n)=\displaystyle\frac1{n\log 2} \log
\left(\sum_{x\in\mathcal{F}_n} \mu(x)^p\right).
 $$
Let $\mu$ and $\nu$ be two probability measures on $\mathcal{C}$
such that, $\nu\ll\mu$ and $\mu$ is a quasi-Bernoulli measure. Then,
$\tau'_\mu(1)$ exists and we have
\begin{eqnarray}\label{E1}
\displaystyle\lim_{ n\to \infty} \frac{\log
\mu\big(B(x,2^{-n})\big)} {\log (2^{-n})}=\displaystyle\lim_{ n\to
\infty} \frac{\log_2 \mu\big(I_n(x)\big)} {-n}= -\tau'_\mu(1) \text{
for}\; \mu\text{-a.e. } x\in \mathcal{C},
\end{eqnarray}
and
\begin{eqnarray}\label{E2}
\displaystyle\lim_{ n\to \infty}\;\; \frac{\log_2
\nu\big(I_n(x)\big)} {-n}= -\tau'_\nu(1)= -\tau'_\mu(1)\quad
\text{for}\; \nu\text{-a.e. } x\in \mathcal{C}.
\end{eqnarray}
For more details about \eqref{E1} and \eqref{E2}, the reader can see
\cite{H3,H2}. We have immediately from \eqref{E1} and \eqref{E2}
that the measure $\nu$ is $(q, \mu)$-unidimensional and
 $$
{\dim}_{\mu}^q(\nu)={\Dim}_{\mu}^q(\nu)=(q-1)\tau'_\mu(1)=(q-1)\tau'_\nu(1).
 $$

\section{Projections results}

In this section, we show that the multifractal Hausdorff and packing
dimensions of a measure $\nu$ are preserved under almost every
orthogonal projection. Casually, we briefly recall some basic
definitions and facts which will be repeatedly used in subsequent
developments. Let $m$ be an integer with $0<m<n$ and $G_{n,m}$ the
Grassmannian manifold of all $m$-dimensional linear subspaces of
$\mathbb{R}^n$. Denote by $\gamma_{n,m}$ the invariant Haar measure
on $G_{n,m}$, such that $\gamma_{n,m}(G_{n,m})=1$. For $V\in
G_{n,m}$, define the projection map $\pi_V: \mathbb{R}^n
\longrightarrow V$ as the usual orthogonal projection onto $V$.
Then, the set $\{\pi_V,\; V \in G_{n,m}\}$ is compact in the space
of all linear maps from $\mathbb{R}^n$ to $\mathbb{R}^m$ and the
identification of $V$ with $\pi_V$ induces a compact topology for
$G_{n,m}$. Also, for a Borel probability measure $\mu$ with compact
support on $\mathbb{R}^n$, denoted by $\supp\mu$, and for $V\in
G_{n,m}$, define the projection $\mu_V$ of $\mu$ onto $V$, by
 $$
\mu_V(A)=\mu(\pi_V^{-1}(A))\quad \forall A\subseteq V.
 $$
Since $\mu$ is compactly supported and $\supp\mu_V=\pi_V(\supp\mu)$
for all $V\in G_{n,m}$, then for any continuous function $f:
V\longrightarrow\mathbb {R}$, we have
 $$
\displaystyle\int_V fd\mu_V=\int f(\pi_V(x))d\mu(x)
 $$
whenever these integrals exist. For more details, see for example
\cite{FM, FN, M1, M2, SB, SS}. The convolution is defined, for
$1\leq m< n$ and $r>0$, by
 $$
\begin{array}{llll}
\overline{\phi}_r^m:\; & \mathbb{R}^n &\longrightarrow & \mathbb{R}
\\ & x & \longmapsto & \gamma_{n,m}\Big\{V\in G_{n,m};\; |\pi_V(x)|\leq
r\Big\}.
\end{array}
 $$
Moreover, define
 $$
\begin{array}{llll}
\phi_r^m: & \mathbb{R}^n & \longrightarrow & \mathbb{R} \\
& x & \longmapsto & \min\Big\{1\:,\: r^m|x|^{-m}\Big\}.
\end{array}
 $$
We have that $\phi_r^m(x)$ is equivalent to $\overline{\phi}_r^m
(x)$ and write $\phi_r^m(x)\asymp \overline{\phi}_r^m(x)$.

 \smallskip \smallskip
For a probability measure $\mu$ and for $V\in G_{n,m}$, we have
\begin{equation*}
\label{r}\mu\ast\phi_r^m(x)\asymp\mu\ast
\overline{\phi}_r^m(x)=\int\mu_V(B(x_V,r))dV
\end{equation*}
and
\begin{equation*}
\mu\ast\phi_r^m(x)=\int\min\Big\{1 \:,\: r^m|x-y|^{-m}\Big\}d\mu(y).
\end{equation*}
So, integrating by parts and converting into spherical coordinates
(see \cite{FN}), we obtain
\begin{equation*}
\mu\ast\phi_r^m(x)=mr^m\int_r^{+\infty}u^{-m-1}\mu(B(x,u))du.
\end{equation*}

 \bigskip
We present the tools, as well as the intermediate results, which
will be used in the proofs of our main results. The following
straightforward estimates concern the behaviour of the convolution
$\mu\ast\phi_r^m(x)$ as $r\to0$.
\begin{lemma}\cite{FN}\label{l1}
Let $1\leq m \leq n$ and $\mu$ be a compactly supported Borel
probability measure on $\R^n$. For all $x\in \mathbb{R}^n$, we have
$$
c r^m\leq\mu\ast\phi_r^m(x)
$$
for all sufficiently small $r$, where $c>0$ independent of $r$.
\end{lemma}
\begin{definition}
Let $E\subseteq \mathbb{R}^n$ and $0 < s <+ \infty$. We say that $E$
is $s$-Ahlfors regular if it is closed and if there exists a Borel
measure $\mu$ on $\mathbb{R}^n$ and a constant $1\leq C_E <
+\infty$, such that $\mu(E)>0$ and
 $$
C_E^{-1} r^s\leq \mu(B(x,r))\leq C_E r^s,\quad\text{for all}\;\;
x\in E\;\;\;\text{and}\;\;\; 0<r\leq 1.
 $$
\end{definition}

\begin{lemma}\label{l2} Let $0<m\leq n$
\begin{enumerate}
\item Let $\mu$ be a compactly supported Borel probability measure on
$\R^n$. Then, for all $x\in \mathbb{R}^n$ and $r>0$,
$$
\mu(B(x,r))\leq\mu\ast\phi_r^m(x).
$$
\item Suppose that $\mu$ is a compactly supported Borel probability measure on $\mathbb{R}^n$ with support
contained in an $s$-Ahlfors regular set for some $0 < s \leq m$.
Then, for all $\varepsilon>0$ and $\mu$-almost all $x$ there is
$c>0$ such that
 $$
c~r^{-\varepsilon}\mu(B(x,r))\geq\mu\ast\phi_r^m(x).
 $$
 for
sufficiently small $r$.
\end{enumerate}
\end{lemma}
\noindent{\bf Proof.} The proof of assertion (1) is exactly the same
as that given in \cite{FN}. The assertion (2) is nothing but Lemma
5.8 of \cite{O}.

 \bigskip
We use the properties of $\mu\ast\phi_r^n(x)$ to have a relationship
between the kernels and the projected measures.
\begin{lemma}\cite{FN}\label{l3} Let $1\leq m \leq n$ and $\mu$ be a compactly supported Borel probability measure on
$\R^n$. We have,
\begin{enumerate}
\item Let $\varepsilon>0$. For all $V\in G_{n,m}$, for $\mu$-almost all
$x$ and for sufficiently small $r$,
 $$
r^\varepsilon \mu\ast\phi_r^m(x)\leq\mu_V(B(x_V,r)).
 $$

\item Let $\varepsilon>0$. For $\gamma_{n,m}$-almost all $V\in
G_{n,m}$, for all $x\in\mathbb{R}^n$ and for sufficiently small $r$,
 $$
r^{-\varepsilon} \mu\ast\phi_r^m(x)\geq\mu_V(B(x_V,r)).
 $$
\end{enumerate}
\end{lemma}

Throughout this section, we consider a compactly supported Borel
probability measure  $\mu$ on $\mathbb{R}^n$ with support contained
in an $s$-Ahlfors regular set for some $0\leq s\leq m < n$ and $\nu$
be a compactly supported Borel probability measure on $\mathbb{R}^n$
such that $\supp\nu \subseteq \supp\mu$ and $\nu\ll\mu$.

\bigskip
We introduce the function $\underline{\alpha}_{\mu,\nu}^{q,m}$ and
$\overline{\alpha}_{\mu,\nu}^{q,m}$, by
 $$
\underline{\alpha}_{\mu,\nu}^{q,m}(x)=\displaystyle\liminf_{ r\to
0}\;\; \frac{\log \nu\ast\phi_r^m(x) -q\log\mu\ast\phi_r^m(x)}{\log
r},
 $$
and
 $$
\overline{\alpha}_{\mu,\nu}^{q,m}(x)=\displaystyle\limsup_{ r\to
0}\;\; \frac{\log \nu\ast\phi_r^m(x)-q\log\mu\ast\phi_r^m(x)}{\log
r}.
 $$
\begin{proposition}\label{prop2}
Let $q\in \mathbb{R}$. We have that for $\nu$-almost all $x$
\begin{enumerate}
\item If $q>0$, then
 $$
\underline{\alpha}_{\mu,\nu}^{q,m}(x)=\underline{\alpha}_{\mu,\nu}^{q}(x).
 $$
\item If $q\leq0$ and  $\underline{\alpha}_{\mu,\nu}^{q}(x)\leq m(1-q)$, then
 $$
\underline{\alpha}_{\mu,\nu}^{q,m}(x)=\underline{\alpha}_{\mu,\nu}^{q}(x).
 $$
\end{enumerate}
\end{proposition}
\noindent{\bf Proof.}
\begin{enumerate}
\item We will prove that for $\nu$-almost all $x$, we have $
\underline{\alpha}_{\mu,\nu}^{q,m}(x)\leq\underline{\alpha}_{\mu,\nu}^{q}(x).
$ The proof of the other inequality is similar.

  \bigskip
By using Lemma \ref{l2}, we have
 $$
\log\nu(B(x,r))\leq\log\nu\ast\phi_r^m(x).
 $$
Since $ \nu$ is absolutely continuous with respect to $\mu$ and
$q>0$, we have that for $\nu$-almost all $x$
 $$
-q\big(\log c-\varepsilon\log r+\log\mu(B(x,r))\big)\leq
-q\log\mu\ast\phi_r^m(x).
 $$
So, for $\nu$-almost all $x$,
 $$
\frac{\log\nu(B(x,r))-q\big(\log c-\varepsilon\log
r+\log\mu(B(x,r))\big)}{\log r}\geq
\frac{\log\nu\ast\phi_r^m(x)-q\log\mu\ast\phi_r^m(x)}{\log r}.
 $$
Finally, Letting $\varepsilon \to 0$, we get
$\underline{\alpha}_{\mu,\nu}^{q}(x)\geq
\underline{\alpha}_{\mu,\nu}^{q,m}(x).$

\item The inequality $\underline{\alpha}_{\mu,\nu}^{q,m}(x)\leq m(1-q)$ follows immediately from
Lemma \ref{l1}.

  \bigskip
 By using Lemma \ref{l2}, we have
 $$
\log\nu(B(x,r))\leq\log\nu\ast\phi_r^m(x).
 $$
Since $q\leq0$, then
 $$
-q\log\mu(B(x,r))\leq-q\log\mu\ast\phi_r^m(x).
 $$
It follows that $\underline{\alpha}_{\mu,\nu}^{q,m}(x)\leq
\underline{\alpha}_{\mu,\nu}^{q}(x)$. The proof for other inequality
is similar to that given for assertion (1).
\end{enumerate}

 \bigskip
The following proposition is a consequence of Lemma \ref{l3}.
\begin{proposition}\label{prop1}
Let $q\in \mathbb{R}$. For $\gamma_{n,m}$-almost all $V\in G_{n,m}$
and $\nu$-almost all $x$, we have
 $$
\underline{\alpha}_{\mu_V,\nu_V}^{q}(x_V)=\underline{\alpha}_{\mu,\nu}^{q,m}(x)
 $$
and
 $$
\overline{\alpha}_{\mu_V,\nu_V}^{q}(x_V)=\overline{\alpha}_{\mu,\nu}^{q,m}(x).
 $$
\end{proposition}

 \bigskip
The following theorem presents general relations between  the
multifractal Hausdorff and the multifractal packing dimension of a
measure $\nu$ and that of its orthogonal projections.
\begin{theorem}Let $q\in \mathbb{R}$.
\begin{enumerate}
\item For $\gamma_{n,m}$-almost all $V\in G_{n,m}$, we have
$$\underline{\Dim}_{\mu_V}^q(\nu_V)=\infess\overline{\alpha}_{\mu,\nu}^{q,m}(x)\quad\text{
and}\quad
\overline{\Dim}_{\mu_V}^q(\nu_V)=\supess\overline{\alpha}_{\mu,\nu}^{q,m}(x)$$
where the essential bounds being related to the measure $\nu$.

\item For $\gamma_{n,m}$-almost all $V\in G_{n,m}$, we have
\begin{enumerate}
\item for $q>0$,
 $$
\begin{array}{ccll}
\underline{\dim}_{\mu_V}^q(\nu_V)  = \underline{\dim}_\mu^q(\nu)
\quad\text{and}\quad\overline{\dim}_{\mu_V}^q(\nu_V)
=\overline{\dim}_\mu^q(\nu).
\end{array}
 $$
\item for $q\leq0$ and $\overline{\dim}_\mu^q(\nu)\leq m(1-q),$
 $$
\underline{\dim}_{\mu_V}^q(\nu_V)=\underline{\dim}_\mu^q(\nu)
\quad\text{and}\quad
\overline{\dim}_{\mu_V}^q(\nu_V)=\overline{\dim}_\mu^q(\nu).
 $$
\end{enumerate}
\end{enumerate}
\end{theorem}
\noindent{\bf Proof.} Follows immediately from Propositions
\ref{prop2} and \ref{prop1} and Corollary \ref{cor1}.
\begin{remark} Due to Proposition 5.10 in \cite{O}, the result is
optimal. If in addition, $q=0$, then the results of Falconer and
O'Neil hold (see \cite{FN}).
\end{remark}

\section*{Acknowledgments}
The author would like to thank the referee for carefully reading of
the manuscript and for various suggestions and improvements.

\end{document}